\documentclass{article}

\usepackage{amssymb, amsmath, amsthm, mathrsfs, bm}
\usepackage{hyperref}

\theoremstyle{plain}
\newtheorem{thm}{Theorem}[section]
\newtheorem{lem}[thm]{Lemma}
\newtheorem{cor}[thm]{Corollary}
\theoremstyle{definition}
\newtheorem{defn}[thm]{Definition}
\newtheorem{rem}[thm]{Remark}
\newtheorem{exmp}[thm]{Example}
\newtheorem*{ack}{Acknowledgements}
\newtheorem*{pf}{Proof}

\title{New proofs of the Assmus-Mattson theorem based on the Terwilliger algebra}
\author{Hajime Tanaka \\ \emph{\small Division of Mathematics, Graduate School of Information Sciences,} \\[-1mm] \emph{\small Tohoku University, Sendai, Japan} \\[-1mm] {\small E-mail: htanaka@ims.is.tohoku.ac.jp} \\ {\small Dedicated to Professor Eiichi Bannai on the occasion of his $60$th birthday}}
\date{}

\begin{document}

\maketitle

\begin{abstract}
We use the Terwilliger algebra to provide a new approach to the Assmus-Mattson theorem.
This approach also includes another proof of the minimum distance bound shown by Martin as well as its dual.
\end{abstract}

\section{Introduction}

The Terwilliger algebra \cite{Terwilliger1992JAC,Terwilliger1993JACa,Terwilliger1993JACb} is an active area of research.
See \cite{Terwilliger2005GC} and the references therein.
The purpose of the present paper is to demonstrate how the theory of the Terwilliger algebra can also be applied to problems in coding theory.

The \emph{Assmus-Mattson theorem} is a very famous theorem relating linear codes and combinatorial designs:
\begin{thm}[Assmus-Mattson \cite{AM1969JCT}]\label{thm:Assmus-Mattson}
Let $Y$ denote a linear code of length $D$ over $\mathbb{F}_q$ with minimum weight $\delta$.
Let $Y^{\perp}$ denote the dual code of $Y$, with minimum weight $\delta^*$.
Suppose $t\in\{1,2,\dots,D\}$ is such that there are at most $\delta-t$  weights of $Y^{\perp}$ in $\{1,2,\dots,D-t\}$, or such that there are at most $\delta^*-t$  weights of $Y$ in $\{1,2,\dots,D-t\}$.
Then the supports of the words of any fixed weight in $Y$ form a $t$-design.
\end{thm}
There are several proofs and strengthenings of this theorem.
See \cite{CDS1991IEEE,CD1993SIAMJDM,Simonis1995LAA,Bachoc1999DCC,Tanabe2001DCC} for instance.
Delsarte \cite{Delsarte1977PRR} proved an Assmus-Mattson-type theorem for general cometric schemes, and Martin \cite{Martin1998JCD} studied the Assmus-Mattson theorem for Johnson schemes based on Delsarte's algebraic version.
Theorem \ref{thm:Assmus-Mattson} has also been generalized to $\mathbb{Z}_4$-linear codes.
See e.g., \cite{Tanabe2003DCC}.

In this paper, we use the Terwilliger algebra to provide a new approach to Theorem \ref{thm:Assmus-Mattson}.
In fact, we prove three versions of the Assmus-Mattson theorem (Theorems \ref{thm:Assmus-Mattson:P}, \ref{thm:Assmus-Mattson:Q}, \ref{thm:Assmus-Mattson:PQ}) and two corollaries (Corollaries \ref{cor:Assmus-Mattson:P}, \ref{cor:Assmus-Mattson:Q}).
Theorem \ref{thm:Assmus-Mattson:Q} coincides with Delsarte's version whereas Theorem \ref{thm:Assmus-Mattson:P} seems new and is the dual to Theorem \ref{thm:Assmus-Mattson:Q} for general metric schemes.
Both theorems are proved by using only the basic properties of the irreducible modules of the Terwilliger algebra.
Corollaries \ref{cor:Assmus-Mattson:P} and \ref{cor:Assmus-Mattson:Q} may improve these theorems assuming sufficient thinness and dual thinness, respectively.
Section \ref{sec:Assmus-Mattson:PQ} deals with metric and cometric schemes.
The main theorem in this section is Theorem \ref{thm:Assmus-Mattson:PQ}, and we apply recent results of Terwilliger on the displacement and split decompositions \cite{Terwilliger2005GC}.
An interesting consequence is that Theorem \ref{thm:Assmus-Mattson} still holds for nonlinear codes as well (with an appropriate interpretation of the ``weights of the dual code''; see Example \ref{exmp:Hamming}).
This is explained in Example \ref{exmp:interpretation}.
In Section \ref{sec:comparisons}, we compare these Assmus-Mattson theorems and their corollaries.
This section also includes a new proof to the minimum distance bound shown by Martin \cite{Martin2000DCC} as well as its dual (Examples \ref{exmp:Martin:P}, \ref{exmp:Martin:Q}).

\section{Preliminaries}

Throughout this paper, let $(X,\bm{R})$ denote a symmetric association scheme with $D$ classes.
Thus $X$ is the vertex set and $\bm{R}=\{R_0,R_1,\dots,R_D\}$ is the set of associate classes.
We refer the reader to \cite{BI1984B,BCN1989B,Terwilliger1992JAC} for terminology and background materials on association schemes and the Terwilliger algebra.

Let $V$ denote a vector space over $\mathbb{C}$ with a distinguished basis $\{\hat{x}:x\in X\}$ and a Hermitian inner product $\langle\hat{x},\hat{y}\rangle=\delta_{xy}$ $(x,y\in X)$.
For every $\chi\in V$ and a subspace $W\subseteq V$, $\chi|_W$ will denote the orthogonal projection of $\chi$ on $W$.
Let $\mathrm{Mat}_X(\mathbb{C})$ denote the $\mathbb{C}$-algebra of all matrices over $\mathbb{C}$ with rows and columns indexed by $X$.
Then $\mathrm{Mat}_X(\mathbb{C})$ acts on $V$ from the left in an obvious manner.
Let $A_0=I,A_1,\dots,A_D\in\mathrm{Mat}_X(\mathbb{C})$ denote the associate matrices and let $E_0=|X|^{-1}J,E_1,\dots,E_D$ denote the primitive idempotents for the Bose-Mesner algebra $M=\langle A_0,A_1,\dots,A_D\rangle$, where $J$ denotes the all ones matrix.

Pick any $x\in X$.
For each $0\leqslant i\leqslant D$, $R_i(x)=\{y\in X:(x,y)\in R_i\}$ will denote the $i$th subconstituent of $(X,\bm{R})$ with respect to $x$.
Let $E_i^*(x), A_i^*(x)\in\mathrm{Mat}_X(\mathbb{C})$ denote the $i$th \emph{dual idempotent} and the $i$th \emph{dual associate matrix with respect to} $x$, respectively ($0\leqslant i\leqslant D$).
They span the \emph{dual Bose-Mesner algebra} $M^*(x)$ \emph{with respect to} $x$.
We recall $E_i^*(x)$, $A_i^*(x)$ are the diagonal matrices with $(y,y)$-entries $(E_i^*(x))_{yy}=(A_i)_{xy}$, $(A_i^*(x))_{yy}=|X|(E_i)_{xy}$.
The \emph{Terwilliger algebra} $T(x)$ \emph{of} $(X,\bm{R})$ \emph{with respect to} $x$ is the subalgebra of $\mathrm{Mat}_X(\mathbb{C})$ generated by $M$ and $M^*(x)$.
We remark $T(x)$ is semisimple and any two nonisomorphic irreducible $T(x)$-modules in $V$ are orthogonal.

For the remainder of this section, fix $x\in X$ and write $E_i^*=E_i^*(x)$ $(0\leqslant i\leqslant D)$, $T=T(x)$.
Let $W\subseteq V$ denote an irreducible $T$-module.
Set
\begin{equation*}
	W_s=\{0\leqslant i\leqslant D:E_i^*W\ne 0\}, \quad W_s^*=\{0\leqslant j\leqslant D:E_jW\ne 0\}.
\end{equation*}
We call $W_s, W_s^*$ the \emph{support} and the \emph{dual support} of $W$, respectively.
The \emph{diameter} (resp. the \emph{dual diameter}) of $W$ is defined by $d(W)=|W_s|-1$ (resp. $d^*(W)=|W_s^*|-1$).
We say $W$ is \emph{thin} whenever $\dim E_i^*W\leqslant 1$ for all $i$, and we say $W$ is \emph{dual thin} whenever $\dim E_jW\leqslant 1$ for all $j$.

Suppose for the moment that $(X,\bm{R})$ is metric with respect to the ordering $A_0,A_1,\dots,A_D$.
The \emph{endpoint} of an irreducible $T$-module $W\subseteq V$ is $r(W)=\min\{0\leqslant i\leqslant D:E_i^*W\ne 0\}$.
We remark that the primary module $M\hat{x}$ is a unique irreducible $T$-module with endpoint zero.
We shall freely use the following basic fact:
\begin{lem}[{\cite[Lemma 3.9]{Terwilliger1992JAC}}]\label{lem:P-property}
Suppose $(X,\bm{R})$ is metric with respect to the ordering $A_0,A_1,\dots,A_D$ and write $A=A_1$.
Let $W\subseteq V$ denote an irreducible $T$-module and set $r=r(W)$, $d=d(W)$.
Then the following hold:
\begin{enumerate}
\item $AE_i^*W\subseteq E_{i-1}^*W+E_i^*W+E_{i+1}^*W$ $(0\leqslant i\leqslant D)$, where $E_{-1}^*=E_{D+1}^*=0$.
\item $W_s=\{r,r+1,\dots,r+d\}$.
\item $E_i^*AE_k^*W\ne 0$ if $|i-k|=1$ $(r\leqslant i,k\leqslant r+d)$.
\item If $W$ is thin, then $W$ is dual thin.
\end{enumerate}
\end{lem}

Next suppose $(X,\bm{R})$ is cometric with respect to the ordering $E_0,E_1,\dots,E_D$.
The \emph{dual endpoint} of an irreducible $T$-module $W\subseteq V$ is $r^*(W)=\min\{0\leqslant j\leqslant D:E_jW\ne 0\}$.
We remark that $M\hat{x}$ is a unique irreducible $T$-module with dual endpoint zero.
\begin{lem}[{\cite[Lemma 3.12]{Terwilliger1992JAC}}]\label{lem:Q-property}
Suppose $(X,\bm{R})$ is cometric with respect to the ordering $E_0,E_1,\dots,E_D$ and write $A^*=A_1^*(x)$.
Let $W\subseteq V$ denote an irreducible $T$-module and set $r^*=r^*(W)$, $d^*=d^*(W)$.
Then the following hold:
\begin{enumerate}
\item $A^*E_jW\subseteq E_{j-1}W+E_jW+E_{j+1}W$ $(0\leqslant j\leqslant D)$, where $E_{-1}=E_{D+1}=0$.
\item $W_s^*=\{r^*,r^*+1,\dots,r^*+d^*\}$.
\item $E_jA^*E_{\ell}W\ne 0$ if $|j-\ell|=1$ $(r^*\leqslant j,\ell\leqslant r^*+d^*)$.
\item If $W$ is dual thin, then $W$ is thin.
\end{enumerate}
\end{lem}

To avoid triviality, we say a vector $\chi\in V$ is a \emph{code} whenever $\chi\not\in E_0V$ and $\chi\not\in E_0^*(z)V$ for every $z\in X$.
We also say a subset $Y\subseteq X$ is a \emph{code} provided its \emph{characteristic vector} $\chi_Y=\sum_{y\in Y}\hat{y}$ is a code; in other words, $Y$ is a code whenever $1<|Y|<|X|$.
To each code $\chi$ in $V$ we associate four fundamental parameters (with respect to the base vertex $x\in X$ and given orderings of the associate matrices and the primitive idempotents):
\begin{gather*}
	\delta_x(\chi)=\min\{i\ne 0:E_i^*\chi\ne 0\},\quad s_x(\chi)=|\{i\ne 0:E_i^*\chi\ne 0\}|, \\
	\delta^*(\chi)=\min\{j\ne 0:E_j\chi\ne 0\},\quad s^*(\chi)=|\{j\ne 0:E_j\chi\ne 0\}|.
\end{gather*}
When $\chi=\chi_Y$ for a code $Y\subseteq X$, we write $\delta_x(Y)$, $s_x(Y)$, and so on.
In this case, we also set
\begin{equation*}
	\delta(Y)=\min\{i\ne 0:\langle\chi,A_i\chi\rangle\ne 0\},\quad s(Y)=|\{i\ne 0:\langle\chi,A_i\chi\rangle\ne 0\}|.
\end{equation*}
We call $\delta(Y)$, $\delta^*(Y)$, $s(Y)$, $s^*(Y)$ the \emph{minimum distance}, \emph{dual distance}, \emph{degree} and \emph{dual degree} of $Y$, respectively.

Some of the most important families of association schemes are associated with \emph{regular semilattices} (see \cite{Delsarte1976JCTA} for the definition).
Below we give two examples:

\begin{exmp}\label{exmp:Hamming}
Let $\mathcal{Q}=\{0,1,\dots,q-1\}$ $(q\geqslant 2)$.
Introduce a new symbol ``$\cdot$'', and let $\mathscr{L}$ denote the set of words of length $D$ over $\mathcal{Q}\cup\{\cdot\}$.
For $u=(u_1,\dots,u_D), v=(v_1,\dots,v_D)\in\mathscr{L}$, we set $u\preccurlyeq v$ if and only if $u_i=\cdot$ or $u_i=v_i$, for all $i$.
Then $(\mathscr{L},\preccurlyeq)$ defines a regular semilattice (\emph{Hamming lattice}) with rank function $\mathrm{rank}(u)=|\{i:u_i\ne\cdot\}|$, and the top fiber induces the Hamming scheme $\mathrm{H}(D,q)$.
It is both metric and cometric.
Every irreducible $T$-module $W\subseteq V$ is thin (thus dual thin) and satisfies $r(W)=r^*(W)$.
Moreover, if $q=2$ then $d(W)(=d^*(W))=D-2r(W)$.
More detailed information on the irreducible $T$-modules of the Hamming scheme can be found in \cite[Section 6]{Terwilliger1993JACb}, \cite{Go2002EJC}.
See also \cite{Schrijver2005IEEE,GST2006JCTA}.
We remark that if $Y$ denotes a code in $\mathrm{H}(D,q)$ then $\delta^*(Y)-1$ coincides with the (maximum) strength of $Y$ as an orthogonal array \cite[Theorem 4.4]{Delsarte1973PRRS}.
If moreover $\mathcal{Q}$ is the finite field $\mathbb{F}_q$, $Y$ is linear with dual code $Y^{\perp}$ and the base vertex $x$ is the zero vector $(0,0,\dots,0)$, then $E_j\chi_Y\ne 0$ if and only if $E_j^*\chi_{Y^{\perp}}\ne 0$ $(0\leqslant j\leqslant D)$ \cite[Chapter 6]{Delsarte1973PRRS}.
See also \cite[Section 2.10]{BCN1989B}.
\end{exmp}

\begin{exmp}
Let $\Omega=\{1,2,\dots,N\}$ and set $\mathscr{L}=\{u\subseteq\Omega:|u|\leqslant D\}$, where $D\leqslant\lfloor N/2\rfloor$.
Then $(\mathscr{L},\preccurlyeq)$, where the partial order $\preccurlyeq$ is given by inclusion, forms a regular semilattice (\emph{truncated Boolean lattice}) with rank function $\mathrm{rank}(u)=|u|$.
The top fiber induces the Johnson scheme $\mathrm{J}(N,D)$.
It is both metric and cometric.
Every irreducible $T$-module $W\subseteq V$ is thin (thus dual thin) and satisfies $r(W)\leqslant r^*(W)$.
Information on the irreducible $T$-modules of the Johnson scheme can be found in \cite[Section 6]{Terwilliger1993JACb}.
See also \cite[Section III]{Schrijver2005IEEE}.
We remark that if $Y$ denotes a code in $\mathrm{J}(N,D)$ then $\delta^*(Y)-1$ coincides with the (maximum) strength of $Y$ as a $t$-$(N,D,\lambda)$ design \cite[Theorem 4.7]{Delsarte1973PRRS}.
\end{exmp}

\section{Assmus-Mattson theorem for metric schemes}

In this section, we assume that $(X,\bm{R})$ is metric with respect to the ordering $A_0,A_1,\dots,A_D$.
Thus $\Gamma=(X,R_1)$ is a distance-regular graph and $\partial(\cdot,\cdot)$ will denote the graph distance in $\Gamma$.
We fix $x\in X$ and write $E_i^*=E_i^*(x)$ $(0\leqslant i\leqslant D)$, $T=T(x)$.
\begin{defn}
For convenience, we say a vector $\chi\in V$ is a \emph{relative} $t$-\emph{codesign with respect to} $x$ if $E_i^*\chi$ and $A_i\hat{x}$ are linearly dependent for all $1\leqslant i\leqslant t$.
\end{defn}
The first version of our Assmus-Mattson theorems is a variant of Delsarte's result (see the remark below):

\begin{thm}[Assmus-Mattson, Version 1]\label{thm:Assmus-Mattson:P}
Let $\chi$ denote a code in $V$.
Set $\delta_x=\delta_x(\chi)$, $s^*=s^*(\chi)$.
Then $A_{\ell}\chi$ is a relative $(\delta_x-s^*)$-codesign with respect to $x$ for $0\leqslant\ell\leqslant D$.
\end{thm}

\begin{pf}
Set $A=A_1$ and $U=(M\hat{x})^{\perp}$ (the orthogonal complement in $V$).
We observe $U$ is the linear span of all irreducible $T$-modules $W\subseteq V$ with $r(W)>0$.
Set $S=\{j\ne 0: E_j\chi\ne 0\}$. 
Then
\begin{equation*}
	\chi|_U\in\left(\sum_{i=\delta_x}^DE_i^*U\right)\cap\left(\sum_{j\in S}E_jU\right).
\end{equation*}
Since $A$ generates $M$ and takes $s^*(=|S|)$ distinct eigenvalues on $\sum_{j\in S}E_jU$, we find $M\chi|_U$ is spanned by $\chi|_U,A\chi|_U,\dots,A^{s^*-1}\chi|_U$ and thus (cf. Lemma \ref{lem:P-property} (i))
\begin{equation*}
	M\chi|_U\subseteq\sum_{i=\delta_x-s^*+1}^DE_i^*U.
\end{equation*}
This proves $E_i^*M\chi\subseteq\mathbb{C}A_i\hat{x}$ for all $1\leqslant i\leqslant\delta_x-s^*$.
In particular, $E_i^*A_{\ell}\chi\in\mathbb{C}A_i\hat{x}$ for $0\leqslant\ell\leqslant D$.
\qed
\end{pf}

\begin{rem}\label{rem:Delsarte:P}
Let $Y$ denote a code in $X$.
Set $\delta=\delta(Y)$, $s^*=s^*(Y)$.
Delsarte \cite[Theorem 5.11]{Delsarte1973PRRS} showed that $Y$ is $(\delta-s^*)$-regular (i.e., $|Y\cap R_{\ell}(z)|$ depends only on $\ell$ and $\partial(z,Y)=\min\{\partial(z,y):y\in Y\}$ whenever $0\leqslant\partial(z,Y)\leqslant\delta-s^*$).
See also \cite[Theorem 11.1.1]{BCN1989B}.
\end{rem}

The following lemma shows that the code $\chi$ in Theorem \ref{thm:Assmus-Mattson:P} exhibits far stronger regularity if irreducible $T$-modules with small endpoints are thin, and validates the term ``Assmus-Mattson'' above:

\begin{lem}\label{lem:characterization:P}
Let $\chi$ denote a vector in $V$.
Then the following are equivalent:
\begin{enumerate}
\item $\chi$ is orthogonal to every irreducible $T$-module $W\subseteq V$ with $1\leqslant r(W)\leqslant t$.
\item $F\chi$ is a relative $t$-codesign with respect to $x$ for any $F\in T$.
In particular, $A_{\ell}\chi$ is a relative $t$-codesign with respect to $x$ for $0\leqslant\ell\leqslant D$.
\end{enumerate}
Suppose every irreducible $T$-module with endpoint at most $t$ is thin.
Then the second part of (ii) implies (i) (and thus (ii)).
\end{lem}

\begin{pf}
With the same notation as in the proof of Theorem \ref{thm:Assmus-Mattson:P}, (ii) is equivalent to $T\chi|_U\subseteq\sum_{i=t+1}^DE_i^*U$;
in other words,
\begin{equation*}
	T\chi|_W\subseteq\sum_{i=t+1}^{r(W)+d(W)}E_i^*W
\end{equation*}
for every irreducible $T$-module $W\subseteq V$ with $r(W)>0$.
Since $T\chi|_W$ equals $0$ or $W$ according to whether $\chi|_W$ is zero or not, the equivalence of (i) and (ii) follows immediately from this comment.

Let $W\subseteq V$ denote an irreducible $T$-module with $1\leqslant r(W)\leqslant t$ and suppose $W$ is thin.
Then $M\chi|_W$ cannot be a subspace of $\sum_{i=t+1}^{r(W)+d(W)}E_i^*W$ unless $\chi|_W=0$ (cf. Lemma \ref{lem:P-property} (iii)).
This completes the proof.
\qed
\end{pf}

As an application of the technique discussed above, we may improve Theorem \ref{thm:Assmus-Mattson:P} assuming sufficient thinness:

\begin{cor}\label{cor:Assmus-Mattson:P}
Let $\chi$ denote a code in $V$.
Set $\delta_x=\delta_x(\chi)$.
Suppose $t\in\{1,2,\dots,D\}$ is such that
\begin{equation*}
	|\{j\in W_s^*:E_j\chi\ne 0\}|\leqslant\delta_x-r(W)
\end{equation*}
for each irreducible $T$-module $W\subseteq V$ with $1\leqslant r(W)\leqslant t$.
If every irreducible $T$-module with endpoint at most $t$ is thin, then $F\chi$ is a relative $t$-codesign with respect to $x$ for any $F\in T$.
\end{cor}

\begin{pf}
Let $W\subseteq V$ denote an irreducible $T$-module with $1\leqslant r(W)\leqslant t$ and set $S=\{j\in W_s^*: E_j\chi\ne 0\}$.
Then as in the proof of Theorem \ref{thm:Assmus-Mattson:P} we find
\begin{equation*}
	M\chi|_W\subseteq\sum_{i=\delta_x-|S|+1}^{r(W)+d(W)}E_i^*W.
\end{equation*}
However since $W$ is thin and $r(W)\leqslant\delta_x-|S|$, this forces $\chi|_W=0$ (cf. Lemma \ref{lem:P-property} (iii)).
Now the result follows from Lemma \ref{lem:characterization:P}.
\qed
\end{pf}

\begin{exmp}
Suppose $(X,\bm{R})$ is the Hamming scheme $\mathrm{H}(D,q)$.
Let $Y$ denote a code in $X$ and assume $\chi_Y$ satisfies the equivalent conditions (i), (ii) in Lemma \ref{lem:characterization:P}.
Then setting $F=A_{\ell}E_k^*$ $(0\leqslant k,\ell\leqslant D)$ we find $|Y\cap R_k(x)\cap R_{\ell}(z)|$ is independent of $z\in R_t(x)$.
In particular, when $x=(0,0,\dots,0)$ the supports of the codewords of fixed weight $k$ form a $t$-design (in $\mathrm{J}(D,k)$).
We remark that for $q=2$, (half of) Theorem \ref{thm:Assmus-Mattson} follows from Corollary \ref{cor:Assmus-Mattson:P}.
\end{exmp}

\begin{exmp}
Again suppose $(X,\bm{R})$ is the Hamming scheme $\mathrm{H}(D,q)$.
When $q$ is a prime power, there are many nonlinear single-error-correcting perfect codes (containing $(0,0,\dots,0)$; see e.g., \cite{PRV2005DCC}).
They have minimum distance three and dual degree one.
Thus Theorem \ref{thm:Assmus-Mattson:P}, together with Lemma \ref{lem:characterization:P}, shows that these codes support $2$-designs.
\end{exmp}

\begin{exmp}\label{exmp:binary-Golay}
The $[24,12,8]$ extended binary Golay code has covering radius four and is self-dual with weight enumerator $x^{24}+759x^{16}y^8+2576x^{12}y^{12}+759x^8y^{16}+y^{24}$ (where $x,y$ are indeterminates).
Thus Corollary \ref{cor:Assmus-Mattson:P} shows that a coset of weight four supports $1$-designs.
On the other hand, it is well-known that the codewords of a fixed weight form a $5$-design.
\end{exmp}

\begin{exmp}
Suppose $(X,\bm{R})$ is the Johnson scheme $\mathrm{J}(N,D)$.
Let $Y$ denote a code in $X$ and assume $\chi_Y$ satisfies the equivalent conditions (i), (ii) in Lemma \ref{lem:characterization:P}.
Then for every $0\leqslant k\leqslant D$, $\{(x-y,y-x):y\in Y\cap R_k(x)\}$ (which is a subset of $\mathrm{J}(D,k)\otimes\mathrm{J}(N-D,k)$) has the following property;
there exists a constant $\lambda$ such that for any $t$-subsets $\xi\subseteq x$ and $\eta\subseteq\Omega-x$, the number of elements $y\in Y\cap R_k(x)$ satisfying $\xi\subseteq x-y$ and $\eta\subseteq y-x$ is exactly $\lambda$.
Such combinatorial objects are (among other things) studied in detail in \cite{Martin1999DCC}.
\end{exmp}

\begin{exmp}\label{exmp:Witt}
The $5$-$(24,8,1)$ large Witt design has block intersection numbers $4$, $2$ and $0$ so that the minimum distance is four, and it can be checked that this design has dual degree two (more precisely, it is a $\{1,2,3,4,5,7\}$-design; see \cite{CDS1991IEEE}).
Thus if the base vertex $x$ is chosen from the design, then by Theorem \ref{thm:Assmus-Mattson:P} we can take $t=2$ in the previous example.
\end{exmp}

\section{Assmus-Mattson theorem for cometric schemes}

In this section, we assume that $(X,\bm{R})$ is cometric with respect to the ordering $E_0,E_1,\dots,E_D$.
We fix $x\in X$ and write $E_i^*=E_i^*(x)$ $(0\leqslant i\leqslant D)$, $M^*=M^*(x)$, $T=T(x)$.
\begin{defn}[\cite{Delsarte1977PRR}]
A vector $\chi\in V$ is said to be a \emph{relative} $t$-\emph{design with respect to} $x$ if $E_j\chi$ and $E_j\hat{x}$ are linearly dependent for all $1\leqslant j\leqslant t$.
\end{defn}
The notion of relative $t$-designs has a geometric interpretation when the scheme is associated with a regular semilattice:

\begin{exmp}\label{exmp:interpretation:Q}
Suppose $(X,\bm{R})$ is induced on the top fiber $X$ of a short regular semilattice $(\mathscr{L},\preccurlyeq)$.
(A semilattice $(\mathscr{L},\preccurlyeq)$ is \emph{short} if $X\wedge X=\mathscr{L}$.)
In this case, Delsarte \cite[Theorem 9.8]{Delsarte1977PRR} showed that $\chi\in V$ is a relative $t$-design with respect to $x$ if and only if for each object $u\in\mathscr{L}$ such that $\mathrm{rank}(u)=t$, $\sum_{y\in X,u\preccurlyeq y}\langle\chi,\hat{y}\rangle$ depends only on $\mathrm{rank}(x\wedge u)$.
\end{exmp}

We may also use the Terwilliger algebra to give a new proof of Delsarte's algebraic version of the Assmus-Mattson theorem (cf. \cite[Theorem 2.8.1]{BCN1989B}):

\begin{thm}[Assmus-Mattson, Version 2 {\cite[Theorem 8.4]{Delsarte1977PRR}}]\label{thm:Assmus-Mattson:Q}
Let $\chi$ denote a code in $V$.
Set $\delta^*=\delta^*(\chi)$, $s_x=s_x(\chi)$.
Then $E_k^*\chi$ is a relative $(\delta^*-s_x)$-design with respect to $x$ for $0\leqslant k\leqslant D$.
\end{thm}

\begin{pf}
By an argument similar to the proof of Theorem \ref{thm:Assmus-Mattson:P}, we find $E_jM^*\chi\subseteq\mathbb{C}E_j\hat{x}$ for all $1\leqslant j\leqslant\delta^*-s_x$.
In particular, $E_jE_k^*\chi\in\mathbb{C}E_j\hat{x}$ for $0\leqslant k\leqslant D$.
\qed
\end{pf}

\begin{rem}\label{rem:Delsarte:Q}
Let $Y$ denote a code in $X$.
Set $\delta^*=\delta^*(Y)$, $s=s(Y)$.
Delsarte \cite[Theorem 5.24]{Delsarte1973PRRS} also showed that if $\delta^*\geqslant s$ then $Y$ is regular (i.e., $0$-regular; for each $k$, $|Y\cap R_k(y)|$ does not depend on the choice of $y\in Y$).
\end{rem}

The following lemma is the counterpart to Lemma \ref{lem:characterization:P}:

\begin{lem}\label{lem:characterization:Q}
Let $\chi$ denote a vector in $V$.
Then the following are equivalent:
\begin{enumerate}
\item $\chi$ is orthogonal to every irreducible $T$-module $W\subseteq V$ with $1\leqslant r^*(W)\leqslant t$.
\item $F\chi$ is a relative $t$-design with respect to $x$ for any $F\in T$.
In particular, $E_k^*\chi$ is a relative $t$-design with respect to $x$ for $0\leqslant k\leqslant D$.
\end{enumerate}
Suppose every irreducible $T$-module with dual endpoint at most $t$ is dual thin.
Then the second part of (ii) implies (i) (and thus (ii)).
\end{lem}

\begin{pf}
Similar to the proof of Lemma \ref{lem:characterization:P}.
\qed
\end{pf}

Assuming sufficient dual-thinness, we may improve Theorem \ref{thm:Assmus-Mattson:Q} as follows:

\begin{cor}\label{cor:Assmus-Mattson:Q}
Let $\chi$ denote a code in $V$.
Set $\delta^*=\delta^*(\chi)$.
Suppose $t\in\{1,2,\dots,D\}$ is such that
\begin{equation*}
	|\{i\in W_s: E_i^*\chi\ne 0\}|\leqslant\delta^*-r^*(W)
\end{equation*}
for each irreducible $T$-module $W\subseteq V$ with $1\leqslant r^*(W)\leqslant t$.
If every irreducible $T$-module with dual endpoint at most $t$ is dual thin, then $F\chi$ is a relative $t$-design with respect to $x$ for any $F\in T$.
\end{cor}

\begin{pf}
Similar to the proof of Corollary \ref{cor:Assmus-Mattson:P}.
\qed
\end{pf}

\begin{exmp}
Suppose $(X,\bm{R})$ is the Hamming scheme $\mathrm{H}(D,q)$.
We remark that in this case the equivalent conditions (i), (ii) in Lemma \ref{lem:characterization:Q} are also equivalent to the conditions (i), (ii) in Lemma \ref{lem:characterization:P}.
For $q=2$, (half of) Theorem \ref{thm:Assmus-Mattson} follows from Corollary \ref{cor:Assmus-Mattson:Q}.
\end{exmp}

\begin{exmp}
Relative $t$-designs as well as the Assmus-Mattson theorem in the Johnson scheme $\mathrm{J}(N,D)$ are studied in detail in \cite{Martin1998JCD} in the context of \emph{mixed block designs}.
(Mixed $2$-designs form a subclass of \emph{balanced bipartite block designs}.)
\end{exmp}

\begin{exmp}
As observed in \cite[Example 10.2]{Delsarte1977PRR} and \cite[Example 9]{Martin1998JCD}, the $5$-$(24,8,1)$ large Witt design provides two relative $3$-designs.
See Example \ref{exmp:Witt}.
\end{exmp}

\section{Assmus-Mattson theorem for metric and cometric schemes}\label{sec:Assmus-Mattson:PQ}

In this section, we assume that $(X,\bm{R})$ is metric with respect to the ordering $A_0,A_1,\dots,A_D$ and cometric with respect to the ordering $E_0,E_1,\dots,E_D$.
We fix $x\in X$ and write $E_i^*=E_i^*(x)$ $(0\leqslant i\leqslant D)$, $T=T(x)$.
The third version of our Assmus-Mattson theorems is related to the displacement and split decompositions for $(X,\bm{R})$ \cite{Terwilliger2005GC}.

Let $W\subseteq V$ denote an irreducible $T$-module.
Then $d(W)=d^*(W)$ by \cite[Corollary 3.3]{Pascasio2002EJC} and moreover \cite[Lemmas 5.1, 7.1]{Caughman1999DM}
\begin{equation*}
	2r(W)+d(W)\geqslant D, \quad 2r^*(W)+d(W)\geqslant D.
\end{equation*}
The \emph{displacement} of $W$ is
\begin{equation*}
	\eta(W)=r(W)+r^*(W)+d(W)-D.
\end{equation*}
Then $0\leqslant \eta(W)\leqslant D$ \cite[Lemma 4.2]{Terwilliger2005GC} and it follows from the above inequalities that $\eta(W)=0$ if and only if $r(W)=r^*(W)=(D-d(W))/2$.
Note that the primary module $M\hat{x}$ has displacement zero.
For each $0\leqslant\eta\leqslant D$, we let $V_{\eta}$ denote the subspace of $V$ spanned by the irreducible $T$-modules with displacement $\eta$.
Then $V=\sum_{\eta=0}^DV_{\eta}$ (orthogonal direct sum) \cite[Lemma 4.4]{Terwilliger2005GC}.
This is the \emph{displacement decomposition of} $V$ \emph{with respect to} $x$.

On the other hand, for $0\leqslant i,j\leqslant D$ we define
\begin{equation*}
	V_{ij}=\left(\sum_{k=0}^iE_k^*V\right)\cap\left(\sum_{\ell=0}^jE_{\ell}V\right).
\end{equation*}
Obviously $V_{i-1,j}, V_{i,j-1}\subseteq V_{ij}$, and we let $\tilde{V}_{ij}$ denote the orthogonal complement of $V_{i,j-1}+V_{i-1,j}$ in $V_{ij}$.
Then for $0\leqslant k,\ell\leqslant D$, we have \cite[Theorem 5.7]{Terwilliger2005GC}
\begin{equation*}
	V_{k\ell}=\sum_{i=0}^k\sum_{j=0}^{\ell}\tilde{V}_{ij} \quad\text{(direct sum)},
\end{equation*}
and in particular, $V=\sum_{i=0}^D\sum_{j=0}^D\tilde{V}_{ij}$ \cite[Corollary 5.8]{Terwilliger2005GC}.
We call the latter sum the \emph{split decomposition of} $V$ \emph{with respect to} $x$.
This decomposition is not orthogonal in general.

Moreover, Terwilliger \cite[Theorem 6.2]{Terwilliger2005GC} showed that for each $0\leqslant\eta\leqslant D$ we have $V_{\eta}=\sum_{i,j}\tilde{V}_{ij}$, where the sum is over $0\leqslant i,j\leqslant D$ such that $i+j=D+\eta$, and thus comparing the displacement and split decompositions he found $\tilde{V}_{ij}=V_{ij}=0$ if $i+j<D$.

\begin{lem}\label{lem:sums}
We have $V_0=\sum_{i=0}^DV_{i,D-i}$ (direct sum).
Moreover the following hold.
\begin{enumerate}
\item $\sum_{k=0}^iV_{k,D-k}=\sum_{k=0}^iE_k^*V_0$ $(0\leqslant i\leqslant D)$.
\item $\sum_{\ell=0}^jV_{D-\ell,\ell}=\sum_{\ell=0}^jE_{\ell}V_0$ $(0\leqslant j\leqslant D)$.
\end{enumerate}
\end{lem}

\begin{pf}
As $V_{i-1,D-i}=V_{i,D-i-1}=0$, we find $V_{i,D-i}=\tilde{V}_{i,D-i}$ and the first line follows.

(i) Clearly it suffices to show $\sum_{k=0}^iE_k^*V_0\subseteq\sum_{k=0}^iV_{k,D-k}$.
Set $A^*=A_1^*(x)$ and let $\theta_k^*$ denote the eigenvalue of $A^*$ associated with $E_k^*V$ $(0\leqslant k\leqslant D)$.
Then it is easy to see that $(A^*-\theta_k^*I)V_{k,D-k}\subseteq V_{k-1,D-k+1}$ for $0<k\leqslant D$ and $(A^*-\theta_0^*I)V_{0D}=0$ (cf. \cite[Theorem 7.1]{Terwilliger2005GC}).
Thus setting $F=\prod_{k=i+1}^D(A^*-\theta_k^*I)$, we obtain
\begin{equation*}
	\sum_{k=0}^iE_k^*V_0=FV_0=\sum_{k=0}^DFV_{k,D-k}\subseteq\sum_{k=0}^iV_{k,D-k}.
\end{equation*}
See also \cite[Theorem 4.6]{ITT2001P}.

(ii) Similar to the proof of (i) above.
\qed
\end{pf}

\begin{thm}[Assmus-Mattson, Version 3]\label{thm:Assmus-Mattson:PQ}
Let $\chi$ denote a code in $V$.
Set $\delta_x=\delta_x(\chi)$, $\delta^*=\delta^*(\chi)$.
Suppose $t\in\{1,2,\dots,D\}$ is such that for every $1\leqslant r\leqslant t$ at least one of the following holds:
\begin{align*}
	|\{r\leqslant j\leqslant D-r:E_j\chi\ne 0\}|&\leqslant\delta_x-r, \\
	|\{r\leqslant i\leqslant D-r:E_i^*\chi\ne 0\}|&\leqslant\delta^*-r.
\end{align*}
If every irreducible $T$-module with displacement zero and endpoint at most $t$ is thin (thus dual thin), then the following hold.
\begin{enumerate}
\item For any $F\in T$, $F\chi$ is orthogonal to $V_{i,D-i}\cap(M\hat{x})^{\perp}$ whenever $1\leqslant i\leqslant t$.
\item For any $F\in T$, $F\chi$ is orthogonal to $V_{D-j,j}\cap(M\hat{x})^{\perp}$ whenever $1\leqslant j\leqslant t$.
\end{enumerate}
\end{thm}

\begin{pf}
First, by the hypothesis we find $\chi|_W=0$ for any irreducible $T$-module $W\subseteq V$ with $\eta(W)=0$ and $1\leqslant r(W)(=r^*(W))\leqslant t$, as in the proofs of Corollaries \ref{cor:Assmus-Mattson:P}, \ref{cor:Assmus-Mattson:Q}.

(i) Set $U_0=V_0\cap(M\hat{x})^{\perp}$.
We observe $U_0$ is the linear span of all irreducible $T$-modules $W\subseteq V$ with $\eta(W)=0$ and $r(W)>0$.
Then, in view of Lemma \ref{lem:sums} it suffices to show $T\chi|_{U_0}\subseteq\sum_{i=t+1}^DE_i^*U_0$, or equivalently $T\chi|_W\subseteq\sum_{i=t+1}^{r(W)+d(W)}E_i^*W$ for every irreducible $T$-module $W\subseteq V$ with $\eta(W)=0$ and $r(W)>0$, but this follows immediately from the above comment.

(ii) Similar to the proof of (i) above.
\qed
\end{pf}

\begin{rem}
The assumption on thinness in Theorem \ref{thm:Assmus-Mattson:PQ} is redundant.
In fact, P. Terwilliger (private communication) pointed out that the irreducible $T$-modules with displacement zero are always thin for any metric and cometric schemes.
This (among other things) will be discussed in a future paper.
\end{rem}

\begin{exmp}\label{exmp:interpretation:PQ}
With the same hypothesis as in Theorem \ref{thm:Assmus-Mattson:PQ}, assume moreover $(X,\bm{R})$ is induced on the top fiber $X$ of a short regular semilattice $(\mathscr{L},\preccurlyeq)$.
For each object $u\in\mathscr{L}$, let $\chi_{\succcurlyeq u}=\sum_{y\in X,u\preccurlyeq y}\hat{y}$ denote the characteristic vector of $\{y\in X:u\preccurlyeq y\}$.
It is a standard fact that $\chi_{\succcurlyeq u}\in V_{D-\mathrm{rank}(u),\mathrm{rank}(u)}$ whenever $u\preccurlyeq x$ \cite{Delsarte1976JCTA}.
We remark that each $A_k\hat{x}$ is obviously a relative $D$-design with respect to $x$ and thus also a relative $t$-design with respect to $x$ (cf. \cite[Corollary 9.9]{Delsarte1977PRR} and the remark that follows it).
Therefore, if $u,v\in\mathscr{L}$ are two objects of rank $t$ such that $u,v\preccurlyeq x$, then in view of the geometric interpretation of relative $t$-designs given in Example \ref{exmp:interpretation:Q}, $\chi_{\succcurlyeq u}-\chi_{\succcurlyeq v}$ is orthogonal to $A_k\hat{x}$ for every $0\leqslant k\leqslant D$;
in other words, $\chi_{\succcurlyeq u}-\chi_{\succcurlyeq v}\in V_{D-t,t}\cap(M\hat{x})^{\perp}$.
Now the second part of Theorem \ref{thm:Assmus-Mattson:PQ} implies that for each $F\in T$, $\langle F\chi,\chi_{\succcurlyeq u}\rangle$ is independent of $u\preccurlyeq x$ with rank $t$.
\end{exmp}

\begin{exmp}\label{exmp:interpretation}
Suppose $(X,\bm{R})$ is the Hamming scheme $\mathrm{H}(D,q)$ and set $x=(0,0,\dots,0)$.
Let $Y$ denote a code in $X$ and set $\chi=\chi_Y$ in the previous example.
Then  we find that (the complements of) the supports of the words of fixed weight $k$ in $Y$ form a $t$-design (in $\mathrm{J}(D,k)\cong\mathrm{J}(D,D-k)$) for every $k$.
In particular, the conclusion of the original Assmus-Mattson theorem (Theorem \ref{thm:Assmus-Mattson}) is also true for nonlinear codes as well.
\end{exmp}

\begin{exmp}\label{exmp:ternary-Golay}
The $[12,6,6]$ extended ternary Golay code has covering radius three and is self-dual with weight enumerator $x^{12}+264x^6y^6+440x^3y^9+24y^{12}$ (where $x,y$ are indeterminates).
Thus Theorem \ref{thm:Assmus-Mattson:PQ} shows that a coset of weight three supports $1$-designs.
On the other hand, it is well-known that the codewords of a fixed weight form a $5$-design.
\end{exmp}

\begin{exmp}
Suppose $(X,\bm{R})$ is the Johnson scheme $\mathrm{J}(N,D)$.
Let $Y$ denote a code in $X$ and set $\chi=\chi_Y$ in Example \ref{exmp:interpretation:PQ}.
Then, in this case we find that the multiset $\{x\cap y:y\in Y\cap R_k(x)\}$ (counting repeats) forms a $t$-design (in $\mathrm{J}(D,D-k)$) for every $k$.
\end{exmp}

\begin{exmp}
The $2$-$(56,12,3)$ design constructed in \cite{BH1980JCTA} has intersection numbers $3,2$ and $0$, and thus Theorem \ref{thm:Assmus-Mattson:PQ} provides two $1$-designs.
\end{exmp}

\section{Comparisons of the Assmus-Mattson theorems}\label{sec:comparisons}

All examples of codes $Y$ known to the author and supporting interesting designs in Hamming and Johnson schemes are regular (in view of Remarks \ref{rem:Delsarte:P}, \ref{rem:Delsarte:Q}), and in general taking the base vertex outside $Y$ seems much less effective than taking the base vertex inside $Y$ (but see Examples \ref{exmp:binary-Golay}, \ref{exmp:ternary-Golay}).
In this section, we compare various Assmus-Mattson theorems and their corollaries focusing on regular codes and base vertices chosen from these codes.

First we assume $(X,\bm{R})$ is metric with respect to the ordering $A_0,A_1,\dots,A_D$ and cometric with respect to the ordering $E_0,E_1,\dots,E_D$.
We also assume that $(X,\bm{R})$ is not an ordinary cycle, and $D\geqslant 3$.

The following result was proved by Ito, Tanabe and Terwilliger in the context of \emph{tridiagonal systems}:
\begin{lem}[{\cite[Lemma 4.5]{ITT2001P}}]\label{lem:ITT}
Fix $x\in X$ and write $E_i^*=E_i^*(x)$ $(0\leqslant i\leqslant D)$.
Let $W\subseteq V$ denote an irreducible $T(x)$-module and set $r=r(W)$, $r^*=r^*(W)$, $d=d(W)(=d^*(W))$.
For $0\leqslant i,j\leqslant d$ define
\begin{equation*}
	W_{ij}=\left(\sum_{k=0}^iE_{r+k}^*W\right)\cap\left(\sum_{\ell=j}^dE_{r^*+\ell}W\right),\quad W_{ij}^*=\left(\sum_{\ell=0}^iE_{r^*+\ell}W\right)\cap\left(\sum_{k=j}^dE_{r+k}^*W\right).
\end{equation*}
Then $W_{ij}=W_{ij}^*=0$ if $i<j$.
\end{lem}
We recall a few standard facts about imprimitivity.
If $(X,\bm{R})$ is imprimitive, then $(X,\bm{R})$ is bipartite or antipodal (or both) \cite[Theorem 4.2.1]{BCN1989B}.
Let $\theta_j$ denote the eigenvalue of $A_1$ associated with $E_jV$ $(0\leqslant j\leqslant D)$.
If $(X,\bm{R})$ is bipartite, then $\theta_j=-\theta_{D-j}$ $(0\leqslant j\leqslant D)$ \cite[Theorem 9.6]{Caughman1998GC}.
If $(X,\bm{R})$ is antipodal, then $(X,\bm{R})$ is an antipodal double cover \cite[Theorem 8.2.4]{BCN1989B}.

Next we recall the cosines.
Let $E\in\{E_0,E_1,\dots,E_D\}$ and let $\theta$ denote the corresponding eigenvalue of $A_1$.
Let $\sigma_0,\sigma_1,\dots,\sigma_D$ denote the real numbers defined by $E=|X|^{-1}m\sum_{i=0}^D\sigma_iA_i$, where $m$ denotes the rank of $E$.
Then $\sigma_0=1$ and for $x,y\in X$ such that $(x,y)\in R_i$ we have $\langle E\hat{x},E\hat{y}\rangle=|X|^{-1}m\sigma_i$.
\begin{lem}[{\cite[Proposition 4.4.7]{BCN1989B}}]\label{lem:cosine}
If $\sigma_i=\pm 1$ for some $i>0$, then one of the following holds:
(i) $\theta=\theta_0$; (ii) $(X,\bm{R})$ is bipartite and $\theta=-\theta_0$; (iii) $(X,\bm{R})$ is antipodal and $i=D$.
\end{lem}
(This lemma is valid for any metric association scheme which is not an ordinary cycle.)
We remark that by the above comments we have $\theta=\theta_D$ if (ii) holds above, and $(X,\bm{R})$ is in fact an antipodal double cover if (iii) holds above.

\begin{thm}\label{thm:Martin:P}
Let $Y$ denote a code in $X$.
Set $\delta^*=\delta^*(Y)$.
Suppose $A_{\ell}\chi_Y$ is a relative $t$-codesign with respect to $x$ for every $0\leqslant\ell\leqslant D$ and $x\in Y$ (where $0\leqslant t\leqslant D$).
If for every $x\in Y$ each irreducible $T(x)$-module $W\subseteq V$ with $1\leqslant r^*(W)\leqslant\delta^*$ satisfies $r(W)\leqslant r^*(W)$,  then one of the following holds:
\begin{enumerate}
\item $(X,\bm{R})$ is bipartite and $Y$ is a bipartite half.
\item $(X,\bm{R})$ is an antipodal double cover, and $Y$ forms an antipodal pair.
\item $\delta^*\geqslant t+1$.
\end{enumerate}
\end{thm}

\begin{pf}
Suppose $t\geqslant\delta^*$.
If $\delta^*=D$ then $t=D$; in other words, $\chi_Y\in M\hat{x}$ for every $x\in Y$.
In this case it is easy to see that we have (i).
Now we assume $\delta^*<D$ and show (ii) holds.

For the moment fix $x\in Y$ and write $E_i^*=E_i^*(x)$ $(0\leqslant i\leqslant D)$, $T=T(x)$.
Then by the hypothesis we have $M\chi_Y|_U\subseteq\sum_{i=t+1}^DE_i^*U$, where $U=(M\hat{x})^{\perp}$.
Equivalently, $M\chi_Y|_W\subseteq\sum_{i=t+1}^{r(W)+d(W)}E_i^*W$ for every irreducible $T$-module $W\subseteq V$ with $r(W)>0$.

Let $W\subseteq V$ denote an irreducible $T$-module with $1\leqslant r^*(W)\leqslant\delta^*$.
Then in particular we have $E_{\delta^*}\chi_Y|_W\in\sum_{i=t+1}^{r(W)+d(W)}E_i^*W$.
But by Lemma \ref{lem:ITT}, this implies $E_{\delta^*}\chi_Y|_W=0$ since $\delta^*-r^*(W)<t+1-r(W)$.
Thus we conclude $E_{\delta^*}\chi_Y\in\mathbb{C}E_{\delta^*}\hat{x}$.

Finally, let $x,y\in Y$.
Then since $E_{\delta^*}\chi_Y\ne 0$ we must have $E_{\delta^*}\hat{x}=\pm E_{\delta^*}\hat{y}$.
Since $0<\delta^*<D$, it follows from Lemma \ref{lem:cosine} that $(X,\bm{R})$ is an antipodal double cover and $\{x,y\}$ is an antipodal pair, as desired.
\qed
\end{pf}

\begin{exmp}\label{exmp:comparison}
Suppose $(X,\bm{R})$ is a Hamming scheme $\mathrm{H}(D,q)$ with $q\geqslant 3$.
Let $Y$ denote a regular code in $X$ and set $\delta=\delta(Y)$, $\delta^*=\delta^*(Y)$, $s^*=s^*(Y)$.
We remark $\delta=\delta_x(Y)$ for every $x\in Y$.
Let $t\in\{1,2,\dots,D\}$ denote an integer satisfying the hypothesis in Corollary \ref{cor:Assmus-Mattson:P} (with $\chi=\chi_Y$ and $x\in Y$).
Then by Theorem \ref{thm:Martin:P} we find $\delta^*\geqslant t+1$.
On the other hand, there exists an irreducible $T$-module with (dual) support $\{t,t+1,\dots, D\}$ (cf. \cite[Section 6]{Terwilliger1993JACb}), so that we must have $s^*\leqslant\delta-t$.
Thus in this case Corollary \ref{cor:Assmus-Mattson:P} gives no improvement on Theorem \ref{thm:Assmus-Mattson:P}.
\end{exmp}

\begin{exmp}\label{exmp:comparison'}
Suppose $(X,\bm{R})$ is the binary Hamming scheme $\mathrm{H}(D,2)$.
Let $Y$ denote a regular code in $X$ which is neither a bipartite half nor a repetition code.
Set $\delta=\delta(Y)$, $\delta^*=\delta^*(Y)$, $s^*=s^*(Y)$, and let $t$ be as in the previous example.
Then Theorem \ref{thm:Martin:P} gives $\delta^*\geqslant t+1$.
Next, set $\delta^{\downarrow *}=\min\{j\ne 0:E_{D-j}\chi_Y\ne 0\}$.
Then since every irreducible $T$-module has displacement zero, we can show $\delta^{\downarrow *}\geqslant t+1$ by an argument similar to the proof of Theorem \ref{thm:Martin:P}.
We omit the details.
We find $s^*\leqslant\delta-t$ if $E_D\chi_Y=0$ and $s^*-1\leqslant\delta-t$ if $E_D\chi_Y\ne 0$.
Thus in view of Remark \ref{rem:improvement} below, in this case Corollary \ref{cor:Assmus-Mattson:P} gives no essential improvement on Theorem \ref{thm:Assmus-Mattson:P}.
For binary linear codes, a similar observation was previously done in \cite[Proposition 3]{Dumer1980MN}.
\end{exmp}

\begin{exmp}[{\cite[Theorem 5]{Martin2000DCC}}]\label{exmp:Martin:P}
Suppose $(X,\bm{R})$ is the Hamming scheme $\mathrm{H}(D,q)$.
Let $Y$ denote a code in $X$, and set $\delta=\delta(Y)$, $\delta^*=\delta^*(Y)$, $s^*=s^*(Y)$.
Then by Theorem \ref{thm:Assmus-Mattson:P}, we can take $t=\delta-s^*$ in Theorem \ref{thm:Martin:P}.
Thus we have $\delta\leqslant\delta^*+s^*-1$ unless $Y$ is a binary repetition code.
(The bipartite halves of $\mathrm{H}(D,2)$ satisfy $\delta=2$, $s^*=1$ and $\delta^*=D$.)
We remark that the proof of this inequality in \cite{Martin2000DCC} uses the classification of perfect codes in Hamming schemes, while we have used the information on the irreducible $T$-modules.
\end{exmp}

\begin{exmp}
Suppose $(X,\bm{R})$ is the Johnson scheme $\mathrm{J}(N,D)$.
Then with the same notation as above, in this case we have $\delta\leqslant\delta^*+s^*-1$ unless $N=2D$ and $Y$ forms a complementary pair.
\end{exmp}

\begin{thm}\label{thm:Martin:Q}
Let $Y$ denote a code in $X$.
Set $\delta=\delta(Y)$.
Suppose $E_k^*(x)\chi_Y$ is a relative $t$-design with respect to $x$ for every $0\leqslant k\leqslant D$ and $x\in Y$ (where $0\leqslant t\leqslant D$).
If for every $x\in Y$ each irreducible $T(x)$-module $W\subseteq V$ with $1\leqslant r(W)\leqslant\delta$ satisfies $r^*(W)\leqslant r(W)$,  then one of the following holds:
\begin{enumerate}
\item $(X,\bm{R})$ is bipartite and $Y$ is a bipartite half.
\item $(X,\bm{R})$ is an antipodal double cover, and $Y$ forms an antipodal pair.
\item $\delta\geqslant t+1$.
\end{enumerate}
\end{thm}

\begin{pf}
Suppose $t\geqslant\delta$.
If $\delta=D$ then $t=D$ and thus $\chi_Y\in M\hat{x}$ for every $x\in Y$.
In this case we have (ii).
Now we assume $\delta<D$ and show (i) holds.

Pick any $x\in Y$ such that $E_{\delta}^*(x)\chi_Y\ne 0$.
Then we find $E_{\delta}^*(x)\chi_Y=A_{\delta}\hat{x}$ by an argument similar to the proof of Theorem \ref{thm:Martin:P}.
Apparently $E_{\delta}^*(y)\chi_Y\ne 0$ for every $y\in R_{\delta}(x)$.
Thus continuing the above argument, we conclude $Y$ contains a connected component of $(X,R_{\delta})$.
Since $0<\delta<D$, it follows that $(X,\bm{R})$ is bipartite and $Y$ is a bipartite half, as desired.
\qed
\end{pf}

\begin{exmp}
Suppose $(X,\bm{R})$ is the Hamming scheme $\mathrm{H}(D,q)$.
Let $Y$ denote a regular code in $X$.
When $q=2$ we assume $Y$ is neither a bipartite half nor a repetition code.
Set $\delta=\delta(Y)$, $\delta^*=\delta^*(Y)$, $s=s(Y)$.
We remark $s=s_x(Y)$ for every $x\in Y$.
Let $t\in\{1,2,\dots,D\}$ denote an integer satisfying the hypothesis in Corollary \ref{cor:Assmus-Mattson:Q} (with $\chi=\chi_Y$ and $x\in Y$).
Then by Theorem \ref{thm:Martin:Q} we find $\delta\geqslant t+1$.
If $q=2$ then we can also show the inequality $\delta^{\downarrow}\geqslant t+1$, where $\delta^{\downarrow}=\min\{i\ne 0:\langle\chi,A_{D-i}\chi\rangle\ne 0\}$.
We omit the details.
As in Examples \ref{exmp:comparison}, \ref{exmp:comparison'} we have $s\leqslant\delta^*-t$ except when $q=2$ and $\langle\chi,A_D\chi\rangle\ne 0$,  in which case we have $s-1\leqslant\delta^*-t$.
In view of Remark \ref{rem:improvement} below, in this case Corollary \ref{cor:Assmus-Mattson:Q} gives no essential improvement on Theorem \ref{thm:Assmus-Mattson:Q}.
\end{exmp}

\begin{exmp}\label{exmp:Martin:Q}
Suppose $(X,\bm{R})$ is the Hamming scheme $\mathrm{H}(D,q)$.
Let $Y$ denote a code in $X$, and set $\delta=\delta(Y)$, $\delta^*=\delta^*(Y)$, $s=s(Y)$.
Then by Theorem \ref{thm:Assmus-Mattson:Q}, we can take $t=\delta^*-s$ in Theorem \ref{thm:Martin:Q}.
Thus we have $\delta^*\leqslant\delta+s-1$ unless $Y$ is a bipartite half of $\mathrm{H}(D,2)$.
(Repetition codes of $\mathrm{H}(D,2)$ satisfy $\delta^*=2$, $s=1$ and $\delta=D$.)
This inequality was first shown by W. J. Martin (private communication; in fact he proved a much stronger inequality).
\end{exmp}

\section{Remarks}

\begin{rem}\label{rem:improvement}
We proved our Assmus-Mattson theorems (Theorems \ref{thm:Assmus-Mattson:P}, \ref{thm:Assmus-Mattson:Q}, \ref{thm:Assmus-Mattson:PQ}) and their corollaries (Corollaries \ref{cor:Assmus-Mattson:P}, \ref{cor:Assmus-Mattson:Q}) by projecting the code $\chi$ to the orthogonal complement $U$ of the primary module $M\hat{x}$.
Thus everything still works, for instance, even if we replace $s_x(\chi)$ by $\tilde{s}_x(\chi)=|\{i\ne 0:E_i^*(x)\chi\not\in\mathbb{C}A_i\hat{x}\}|$ and/or $s^*(\chi)$ by $\tilde{s}^*(\chi)=|\{j\ne 0:E_j\chi\not\in\mathbb{C}E_j\hat{x}\}|$.
This (slight) improvement seems particularly effective for codes in the binary Hamming scheme $\mathrm{H}(D,2)$ (in which case $\dim E_D^*(x)V=\dim E_DV=1$).
See also \cite[Section 2.8]{BCN1989B}.
\end{rem}

\begin{rem}
Lalaude-Labayle \cite{Lalaude-Labayle2001IEEE} classified the self-orthogonal binary linear codes with minimum weight at most $10$ (resp. $18$) and whose words of minimum weight support $3$-designs (resp. $5$-designs).
\end{rem}

\begin{rem}
Recently, Schrijver \cite{Schrijver2005IEEE} established the \emph{semidefinite programming bound} on the sizes of codes in the binary Hamming schemes and Johnson schemes, which is shown to be at least as good as Delsarte's bound based on the linear programming method \cite{Delsarte1973PRRS}.
This provides a remarkable application of the Terwilliger algebra.
See also \cite{GST2006JCTA}.
\end{rem}

\begin{ack}
The author would like to thank Jack Koolen, Bill Martin, Akihiro Munemasa and Paul Terwilliger for helpful discussions and comments.
In particular, Paul Terwilliger drastically simplified the initial proof of Corollaries \ref{cor:Assmus-Mattson:P} and \ref{cor:Assmus-Mattson:Q}, which ultimately led to the whole results presented in this paper.
Part of this work was done while the author was visiting the Combinatorial \& Computational Mathematics Center at Pohang University of Science and Technology.
The author wishes to thank the Center for its hospitality during this visit.
The author's research is supported by the Japan Society for the Promotion of Science.
\end{ack}



\small
\begin{thebibliography}{99}

\bibitem{AM1969JCT}
E. F. Assmus, Jr. and H. F. Mattson, Jr.,
New $5$-designs,
J. Combin. Theory 6 (1969) 122-151.

\bibitem{Bachoc1999DCC}
C. Bachoc,
On harmonic weight enumerators of binary codes,
Des. Codes Cryptogr. 18 (1999) 11-28.

\bibitem{BI1984B}
E. Bannai and T. Ito,
Algebraic combinatorics I,
Benjamin/Cummings, Menlo Park, 1984.

\bibitem{BH1980JCTA}
H. Beker and W. Haemers,
$2$-designs having an intersection number $k-n$,
J. Combin. Theory Ser. A 28 (1980) 64-81.

\bibitem{BCN1989B}
A. E. Brouwer, A. M. Cohen and A. Neumaier,
Distance-regular graphs,
Springer-Verlag, Berlin, 1989.

\bibitem{CD1993SIAMJDM}
A. R. Calderbank and P. Delsarte,
On error-correcting codes and invariant linear forms,
SIAM J. Discrete Math. 6 (1993) 1-23.

\bibitem{CDS1991IEEE}
A. R. Calderbank, P. Delsarte and N. J. A. Sloane,
A strengthening of the Assmus-Mattson theorem,
IEEE Trans. Inform. Theory 37 (1991) 1261-1268.

\bibitem{Caughman1998GC}
J. S. Caughman, IV,
Spectra of bipartite $P$- and $Q$-polynomial association schemes,
Graphs Combin. 14 (1998) 321-343.

\bibitem{Caughman1999DM}
J. S. Caughman, IV,
The Terwilliger algebras of bipartite $P$- and $Q$-polynomial schemes,
Discrete Math. 196 (1999) 65-95.

\bibitem{Delsarte1973PRRS}
P. Delsarte,
An algebraic approach to the association schemes of coding theory, Philips Res. Rep. Suppl. No. 10 (1973).

\bibitem{Delsarte1976JCTA}
P. Delsarte,
Association schemes and $t$-designs in regular semilattices,
J. Combinatorial Theory Ser. A 20 (1976) 230-243.

\bibitem{Delsarte1977PRR}
P. Delsarte,
Pairs of vectors in the space of an association scheme,
Philips Res. Rep. 32 (1977) 373-411.

\bibitem{Dumer1980MN}
I. I. Dumer,
Remarks on tactical configurations in codes,
Math. Notes 28 (1980) 543-547.

\bibitem{GST2006JCTA}
D. Gijswijt, A. Schrijver and H. Tanaka,
New upper bounds for nonbinary codes based on the Terwilliger algebra and semidefinite programming,
J. Combin. Theory Ser. A 113 (2006) 1719-1731.

\bibitem{Go2002EJC}
J. T. Go,
The Terwilliger algebra of the hypercube,
European J. Combin. 23 (2002) 399-429.

\bibitem{ITT2001P}
T. Ito, K. Tanabe and P. Terwilliger,
Some algebra related to $P$- and $Q$-polynomial association schemes,
Codes and association schemes (Piscataway, NJ, 1999),
Amer. Math. Soc., Providence, RI (2001), pp. 167-192;
arXiv:\href{http://arxiv.org/abs/math.CO/0406556}{math.CO/0406556}.

\bibitem{Lalaude-Labayle2001IEEE}
M. Lalaude-Labayle,
On binary linear codes supporting $t$-designs,
IEEE Trans. Inform. Theory 47 (2001) 2249-2255.

\bibitem{Martin1998JCD}
W. J. Martin,
Mixed block designs,
J. Combin. Des. 6 (1998) 151-163.

\bibitem{Martin1999DCC}
W. J. Martin,
Designs in product association schemes,
Des. Codes Cryptogr. 16 (1999) 271-289.

\bibitem{Martin2000DCC}
W. J. Martin,
Minimum distance bounds for $s$-regular codes,
Des. Codes Cryptogr. 21 (2000) 181-187.

\bibitem{Pascasio2002EJC}
A. A. Pascasio,
On the multiplicities of the primitive idempotents of a $Q$-polynomial distance-regular graph,
European J. Combin. 23 (2002) 1073-1078.

\bibitem{PRV2005DCC}
K. T. Phelps, J. Rif\`{a} and M. Villanueva,
Kernels and $p$-kernels of $p^r$-ary $1$-perfect codes,
Des. Codes Cryptogr. 37 (2005) 243-261.

\bibitem{Schrijver2005IEEE}
A. Schrijver,
New code upper bounds from the Terwilliger algebra and semidefinite programming,
IEEE Trans. Inform. Theory 51 (2005) 2859-2866.

\bibitem{Simonis1995LAA}
J. Simonis,
MacWilliams identities and coordinate partitions,
Linear Algebra Appl. 216 (1995) 81-91.

\bibitem{Tanabe2001DCC}
K. Tanabe,
A new proof of the Assmus-Mattson theorem for non-binary codes,
Des. Codes Cryptogr. 22 (2001) 149-155.

\bibitem{Tanabe2003DCC}
K. Tanabe,
A criterion for designs in $\mathbb{Z}_4$-codes on the symmetrized weight enumerator,
Des. Codes Cryptogr. 30 (2003) 169-185.

\bibitem{Terwilliger1992JAC}
P. Terwilliger,
The subconstituent algebra of an association scheme I,
J. Algebraic Combin. 1 (1992) 363-388.

\bibitem{Terwilliger1993JACa}
P. Terwilliger,
The subconstituent algebra of an association scheme II,
J. Algebraic Combin. 2 (1993) 73-103.

\bibitem{Terwilliger1993JACb}
P. Terwilliger,
The subconstituent algebra of an association scheme III,
J. Algebraic Combin. 2 (1993) 177-210.

\bibitem{Terwilliger2005GC}
P. Terwilliger,
The displacement and split decompositions for a $Q$-polynomial distance-regular graph,
Graphs Combin. 21 (2005) 263-276;
arXiv:\href{http://arxiv.org/abs/math.CO/0306142}{math.CO/0306142}.

\end{thebibliography}
\end{document}